\documentclass[twocolumn,superscriptaddress]{revtex4-1}

\usepackage{graphicx}
\usepackage{float}
\usepackage{epstopdf}
\usepackage{amsmath}
\usepackage{amssymb}
\usepackage{mdframed}
\usepackage[utf8]{inputenc}
\usepackage[dvipsnames]{xcolor}
\usepackage{tikz}

\usetikzlibrary{automata,positioning,arrows.meta}

\usepackage[caption = false]{subfig}

\usetikzlibrary{automata,positioning,arrows.meta}

\usepackage[caption = false]{subfig}

\makeatletter

\newcommand{\Rmnum}[1]{\expandafter\@slowromancap\romannumeral #1@}

\makeatother

\begin{document}

\title{{Optima and Simplicity in Nature}}

\author{Kamaludin Dingle}
\email{Correspondence: dinglek@caltech.edu}
\affiliation{Department of Chemical Engineering and Biotechnology, Cambridge University, UK}
\affiliation{Department of Computing and Mathematical Sciences, California Institute of Technology, USA}
\affiliation{CAMB,
Department of Mathematics and Natural Sciences, 
Gulf University for Science and Technology, Kuwait}

\date{\today}

\begin{abstract}
\noindent
Why are simple, regular, and symmetric shapes common in nature? Many natural shapes arise as solutions to energy minimisation or other optimisation problems, but is there a general relation between optima and simple, regular shapes and geometries? Here we argue from algorithmic information theory that for objective functions common in nature ---  based on physics and engineering laws --- optimal geometries will be simple, regular, and symmetric. Further, we derive a null model prediction that if a given geometry is an optimal solution for one natural objective function, then it is \emph{a priori} more likely to be optimal or close to optimal for another objective function. 
\end{abstract}


\maketitle
\twocolumngrid

\section{Introduction}
Certain simple shapes, such as spheres, spirals, and fractal branching patterns, abound in nature. Why are they so common? One perspective on geometries in nature is that they arise via solutions to some form of optimisation problem, where the objective function stems from basic physics and engineering principles \cite{hildebrandt1996parsimonious}. This suggests there might be some general relation between optima and simplicity \cite{cohn2010order}; but what precisely is the connection, if any? There are also several examples of simple, regular, or symmetric geometries that are simultaneously optimal, or near optimal, for different objective functions. Consider these motivating examples: \\
(a) There exists a universally optimal arrangements of points on a sphere, which is optimal for a  range of different energy functions \cite{cohn2007universally};\\
(b) In studies of lattice proteins --- a simplified biophysics model of protein folding --- it has been observed that simple and regular lattice protein folds are the most designable (ie have the most sequences assigned) \cite{li1996emergence}, are the fastest folders \cite{melin1999designability}, the fastest unfolders \cite{dias2006designable}, the most thermally stable \cite{li1996emergence}, and are highly robust to genetic mutations \cite{greenbury2016genetic} (see also \cite{nelson1997symmetry});\\
(c) Scale free network architectures, with their regular patterns have been found to be highly robust \cite{valverde2002scale}, and also efficient \cite{cohen2003scale} due to having small network diameters; \\
(d) Spheres represent the optimal convex 3D shape for reduced thermal loss due to minimising surface area for a given volume, but at the same time this shape is believed to be optimal for low density packing of convex shapes \cite{haji2009disordered}. 


Here we adopt an algorithmic information theory \cite{solomonoff1960preliminary,kolmogorov1965three,chaitin1975theory} approach, and argue that if the objective function to be optimised is `simple' (in a technical sense, discusses below) and the set of geometries over which the the objective function is evaluated is also `simple', then:\\
\noindent
(1) \emph{Optimal geometries will be simple, symmetric, or regular}
Further, perhaps more surprisingly we derive a null-model expectation that: \\
\noindent
(2) \emph{Given that a geometry is optimal for one function, then it is more likely to also be an optimal geometry for another function, compared to a null model based on random functions}\\
We also argue that these types of optimisation problem occur commonly in nature, hence explaining at least some of the ubiquitous appearances of simple natural shapes and forms. Before deriving these statements, in the next section we give some brief theoretical background on algorithmic information theory and formally describe the the optimisation problems that we will study.

\section{Background and problem set up}

\noindent
{\bf Theoretical framework:} Developed within theoretical computer science, \emph{algorithmic information theory} \cite{solomonoff1960preliminary,kolmogorov1965three,chaitin1975theory} (AIT) sits at the intersection of computation, computability theory, and information theory. AIT is fundamentally based on  estimating the information content, or complexity, of discrete objects or patterns, such as discrete sequences, or geometries. These information estimates can yield mathematical relations and bounds, which can be used for various scientific and mathematical predictions. In AIT, information content is quantified by \emph{Kolmogorov complexity}, $K(x)$, which measures the complexity of an individual object $x$ in terms of the amount of information required to describe or generate $x$. Another way to think about $K(x)$ is via compression, where $K(x)$ measures the size of the compressed version of $x$.  If $x$ is a simple string eg containing repeating patterns like $x=101010101010101010$ then it is easy to compress, and hence $K(x)$ will be small. Contrast this with a randomly generated bit string of length $n$, which is highly unlikely to contain any patterns, and hence cannot be compressed. In general, a sequence of length $n$ is called ``complex'' or ``random'' if $K(x)\approx n$ bits, and ``simple'' if $K(x) \ll n$. 

More formally, the Kolmogorov complexity $K_U(x)$ of a string $x$ with respect to a (prefix optimal) universal Turing machine \cite{turing1936computable} (UTM) $U$,  is defined \cite{solomonoff1960preliminary,kolmogorov1965three,chaitin1975theory} as
\begin{equation}
K_U(x) = \min_{p}\{|p|: U(p)=x\}
\end{equation}
where $p$ is a binary program for $U$, and $|p|$ denotes the length of the (halting) program $p$ in bits.  Due to the invariance theorem \cite{li2008introduction} for any two optimal UTMs $U$ and $V$, $| K_U(x) - K_V(x)| \leq c$ so that the complexity of $x$ is independent of the choice of the machine, to within an additive constant $c$. Hence  we conventionally drop the subscript $U$ in $K_U(x)$, and speak of `the' Kolmogorov complexity $K(x)$. $K(x)$ is formally uncomputable, meaning that there  cannot exist a general algorithm that for any arbitrary string returns the value of $K(x)$. In practice, it is commonly approximated by standard data compression algorithms \cite{li2008introduction}. 
Further, a large number of studies have shown that AIT and Kolmogorov complexity can be successfully applied in physics, including thermodynamics \cite{bennett1982thermodynamics,kolchinsky2020thermodynamic,zurek1989algorithmic}, entropy estimation \cite{avinery2019universal,martiniani2019quantifying},  in addition to applications in engineering and statistics \cite{vitanyi2013similarity,li2008introduction,dingle2018input,dingle2020generic}. More details and technicalities can be found in standard AIT references \cite{li2008introduction,calude2002information,gacs1988lecture,shen2022kolmogorov}.\\  


 \noindent
{\bf Problem description:} We consider a general optimization problem to take the form: find the optimal geometry $x^*$ where
\begin{equation}\label{eq:max}
x^*=\arg\left(\underset{x}\max \left(f(x)\right)\right)
\end{equation}
with $x\in \mathcal{X}$, and $f(x)$ is the objective function  to maximise. 
The set $\mathcal{X}$ is assumed to be some finite set of discrete geometries, configurations, or sequences, such as binary strings of length $n$ bits, graphs on $n$ nodes, the set of self avoiding walks of length $n$ on a lattice, etc. The value of $n$ is a `natural' measure of the size of the problem. We assume that we can enumerate all possible elements within $\mathcal{X}$. Computational complexity is irrelevant for our purposes, so it does not matter if enumerating all geometries would require exponentially large amounts of time. We assume that $f(x)$ is some computable function that can be evaluated eg on a computer. Combinatorial optimization problems often come in this form \cite{mezard2009information}.

\section{Results} 
 
 \noindent
{\bf Bounding the complexity of optima:} We can describe any $x\in\mathcal{X}$, given the value of $n$, by first generating the set $\mathcal{X}$ using $n$, then producing an ordering of the elements of $\mathcal{X}$ using the function $f$ (ie ordering elements from most optimal to least), and then describing $x$ exactly with its rank $r$ in the ordering. Thus we can bound the Kolmogorov complexity of optima by implementing the following algorithm: 
\begin{eqnarray}
&&\text{i) Enumerate in order all elements of $\mathcal{X}$}\nonumber\\
&&\text{ii) For each $x_i\in \mathcal{X}$, evaluate $f(x_i)$}\nonumber\\
&&\text{iii) List $x_i$ in descending order of optimality: $x^{*}_1$, $x^{*}_2$,\dots,$x^{*}_{|\mathcal{X}|}$}\nonumber\\
&&\text{iv) Set $x^*=x^{*}_1$, print $x^*$ and halt}\nonumber
\end{eqnarray}
Any element in the list  $x^{*}_1$, $x^{*}_2$,\dots,$x^{*}_{|\mathcal{X}|}$ can be identified by its rank $r$ in the list, where eg $r=1$ corresponds the optimum $x^*$, $r=2$ indicates the second highest value for $f$, and $r=|\mathcal{X}|$ is assigned to the geometry with the lowest value of the objective function. If evaluating the objective function for different elements leads to the same function value, then this could lead to multiple elements with the same rank $r$. To avoid this and keep all ranks unique, we will assume that elements with the same value are also ranked according to the order they appear in the enumeration. For example, if $f(x_{45})= f(x_{52})$ and $x_{45}$ has rank $r$, then $x_{52}$ will have rank $r+1$ because it appeared later in the original enumeration.

Putting  all the preceding together we can bound the complexity of any geometry $x^{*}_r$ by
\begin{equation}
K(x^{*}_r|n)\leq K(\mathcal{X}|n) + K(f) +K(r|n)+O(1)\label{eq:mainbound}
\end{equation}
From this bound we see that the complexity $K(\mathcal{X}|n)$ of the set $\mathcal{X}$ as well as the complexity $K(f)$ of the function $f$ are important quantities. 
We will now examine these quantities in more detail.

 The set $\mathcal{X}$ may in general be simple or complex. Many optimization problems in physics take the form looking for an optimal $x$ over \emph{all possible} $x$ within some set, such as all possible graphs up to some size, all possible protein configurations, all possible branching patterns, all possible topologies with some condition. Perhaps surprisingly, the ``set of all'' objects of some kind typically has very low Kolmogorov complexity. For example, the set of all binary strings of length $n$ has complexity only $K(n)\lesssim \log_2(n)$, because it can be described by the program  ``Enumerate and print all binary strings of length $n$''.  (Contrast this with the complexity of a typical element of the set of all binary strings of length $n$, which is likely to have complexity $\approx n$ bits, far higher than $\log_2(n)$ bits.) Indeed, the countably infinite set of \emph{all} binary strings \{0,1,00,01,...\} has very low information content, because it can be generated by a fixed length program, and hence has complexity $O(1)$ bits. Similarly, the set of all graphs on $n$ nodes has very low complexity, at most $\approx\log_2(n)$ bits.  A set can be complex also, if many bits are required to specify it precisely. For example if $\mathcal{X}$ is made from randomly selecting $M\ll 2^n$ strings out of all binary strings of length $n$, then in this case $K(\mathcal{X})\approx Mn$ bits because the only way to precisely generate this set is by individually specifying each element, each of which has complexity $\approx n$ bits (with high probability).

An objective function $f(x)$  may in general be simple or complex. Examples of simple functions are (low order) polynomials, trigonometric functions, and many of the equations describing the laws of physics (eg $F=ma$). These all are classed as simple $O(1)$ complexity functions, because they can be computed via fixed sized algorithms.  An example of a complex function is $f(x)=x+0110....11$ where the latter binary string is a randomly chosen string of length $n$. In this case, $K(f)\approx n$ bits, which may be arbitrarily large. 
A noteworthy and less trivial example of a complex function arrises in the well-know \emph{assignment problem} from combinatorial optimisation \cite{mezard2009information}: If there are $n$ people and $n$ jobs, and a matrix $C_{ij}$ which gives the affinity of person $i$ for job $j$, the problem is to find the assignment of the jobs to the people which maximises the total affinity. The space of configurations/sequences if the space of the  $n!$ permutations of the $n$ indices. Here, the objective function $f(x)$ is in general not simple, because the Kolmogorov complexity of the matrix $C_{ij}$ will typically be $O(n^2)$. \\

 \noindent
{\bf Simple optimal geometries:} We have seen that for many optimisation problems relevant to physics, simple sets with simple objective functions can be expected. In these cases, both  $K(\mathcal{X}|n)$ and $K(f)$ will be small. Further, highly optimal structures $x$ will have rank values $r\approx 1$ (and $r=1$ for the optimum) and hence $K(r|n)\approx 0$. Therefore the right hand side of Eq.\ (\ref{eq:mainbound}) will be very small, and 
\begin{equation}
K(x^{*}_r|n)\leq K(\mathcal{X}|n) + K(f) +K(r|n)\approx 0 \Rightarrow K(x^{*}_r|n)\approx 0 \label{eq:simpleoptima}
\end{equation}
so that any $x^{*}_r$ with $r\approx 1$ must be simple, and $x^{*}$ must have very low, $O(1)$ complexity. 

We have proven that highly optimal solutions must have low Kolmogorov complexity, but this is not necessarily the same as simple or symmetric in the usual senses of the words. It is possible, in principle, that these solutions appear complex or pseudo random. However, we can expect in many natural systems that near optimal solutions will in fact be simple, or regular, or symmetric, because common physics-based objective functions are unlikely to contain UTMs which can compute arbitrary algorithms and generate pseudo random patterns. Notable exceptions include chaotic systems, which are common in nature, and hence may possible lead to pseudo random optima. See Vitanyi \cite{vitanyi2013similarity,vitanyi2020incomputable} for a similar argument for why common patterns in natural science are unlikely to be pseudo random.

Note that $K(r|n)\approx 0$ also holds for $r\approx |\mathcal{X}|$. For example if there are $2^n$ possible geometries, then the last rank will be $r=2^n$, and hence $K(r|n)=K(2^n|n)=O(1)$. This also implies the \emph{least} optimal geometry must also have low Kolmogorov complexity and hence be simple (assuming a simple set and function). In general, any geometries with extreme values with respect to the objective function will be simple.\\

\noindent
{\bf Example of a simple optimum:}  Define $\mathcal{X}$ to be all binary strings of length $n$, and let the objective function be $f(x)=\sum^n x_i$, with $x_i$ the $i$th bit of the string $x$ (i.e.\ $f(x)$  counts the number of 1's in the string).  The set $\mathcal{X}$ has complexity $K(\mathcal{X}|n)=O(1)$, and $K(f)=O(1)$. It follows then from Eq.\ (\ref{eq:simpleoptima})  that whatever is the optimal string $x^\star$ with $r=1$, it must have $K(x^\star|n)=O(1)$ and hence be very simple. In this somewhat trivial example, we can confirm this prediction by noting that the optimal solution is $x^\star=111\dots11$ with $n$ 1's, and we know that $K(111\dots11|n)=O(1)$. \\

\noindent
{\bf Null expectation for simultaneously optimising multiple objective functions:} In the Introduction we saw several examples of optimal geometries that simultaneously have high or optimal values for different objective functions.  But are these rare exceptional cases, or  instances of a more common trend? We now derive a surprising null expectation connecting optima for different objective functions.
 
Consider a finite and simple set $\mathcal{X}$ as above, and two different objective functions $f_1$ and $f_2$, with respective optima $x^{*}_{f_1}, x^{*}_{f_2}\in\mathcal{X}$. If these two functions are random functions, such that their values are random variables, then the null model prediction for the probability that the two optimal geometries coincide is
\begin{equation}
P( x^{*}_{f_1}= x^{*}_{f_2}) = \frac{1}{|\mathcal{X}|}
\end{equation}
assuming that all function values are distinct. This probability is typically exponentially small, because $|\mathcal{X}|$ is typically exponential in $n$, eg $|\mathcal{X}|=2^n$ if $\mathcal{X}$ consists of all binary strings. Hence, under the random function null model it is highly unlikely that the optimum geometry for one function coincides with that of another function.

The situation is very different for optima from simple functions: if $f_1$ and $f_2$ are simple having $O(1)$ complexity, then both $x^{*}_{f_1}$ and  $x^{*}_{f_2}$ must be simple by Eq.\ (\ref{eq:simpleoptima}). Because there are exponentially few simple geometries, $x^{*}_{f_1}$ and  $x^{*}_{f_2}$ must both come from this exponentially small set. Hence given that we already know $x^{*}_{f_1}$ is simple and therefore  in this small set, it is \emph{a priori} far more likely that $x^{*}_{f_1}$ and  $x^{*}_{f_2}$ coincide. More formally, let $\mathcal{X}_{O(1)} \subset\mathcal{X}$ be such that all elements are simple with $O(1)$ complexity. In this scenario, the null probability expectation that the two optima coincide is
\begin{equation}
P( x^{*}_{f_1}= x^{*}_{f_2}) = \frac{1}{|\mathcal{X}_{O(1)}|} \gg  \frac{1}{|\mathcal{X}|}
\end{equation}
because $ |\mathcal{X}_{O(1)}|=O(1) \ll |\mathcal{X}|$ due to the fact that there are only $O(1)$ simple geometries. While $O(1)$ is not a precise quantitative value, nonetheless we understand that this probability is much higher than the exponentially small value of $ 1/|\mathcal{X}|$. $1/|\mathcal{X}_{O(1)}$ may not be a high value, but it is much higher than expected by random functions.

The preceding is an \emph{a priori} argument, or expected null model for simple optima. It does not mean that all simple optima will be optima for all simple functions. Indeed,  if we chose  $f_2(x)=-f_1(x)$ then the optimum for $f_1$ will be the \emph{worst} performing geometry for $f_2$. This illustrates that extrema must be simple, but not necessarily optimal.

Another reason why optima may not coincide is due the fact that while optima (under the described conditions) must be simple, it does not mean that all simple geometries must be optima. Some simple geometries may have low or medium values with respect to the objective function. This is related to the occurrence of low complexity, low probability patterns in studies of simplicity bias \cite{dingle2020generic,alaskandarani2022low}.

The above argument can be extended beyond exactly coinciding optima, to geometries which are close to optima ie with rank $r\approx 1$. 
We can ask, given we have an optimal solution $x^{*}_{f_1}$ for $f_1$, and we have a set $\mathcal{X}_{f_2}^{r}\subset\mathcal{X}$ consisting of the $r$ most optimal geometries in $\mathcal{X}$ ranked according to $f_2$, what is the probability that  $x^{*}_{f_1} \in \mathcal{X}_{f_2}^{r}$? This is an interesting question, because if $x^{*}_{f_1} \in \mathcal{X}_{f_2}^{r}$ then it means that the optimal geometry for one objective function $f_1$ is also close to optimal for another apparently unrelated objective function $f_2$. To proceed, if $f_1$ and $f_2$ are random functions, then we expect
\begin{equation}
P(x^{*}_{f_1} \in \mathcal{X}_{f_2}^{r})= \frac{r}{|\mathcal{X}|} 
\end{equation}
which is a very small number for highly optimal geometries with $r\approx 1$ because $r\ll |\mathcal{X}|$. Contrast this with the case where functions are simple, then the null model prediction is 
\begin{equation}
P(x^{*}_{f_1} \in \mathcal{X}_{f_2}^{r})= \frac{r}{|\mathcal{X}_{O(1)}|}\gg \frac{r}{|\mathcal{X}|}
\end{equation}
because both $x^{*}_{f_1}\in \mathcal{X}_{O(1)}$ and $\mathcal{X}_{f_2}^{r}\subset \mathcal{X}_{O(1)}$ due to being simple.

\section{Discussion}
We have shown that optimal discrete geometries (or configurations) selected by low Kolmogorov complexity (simple) objective functions from low Kolmogorov complexity sets, must themselves have low Kolmogorov complexity: optima for simple objective functions must be simple. However, geometries with very \emph{low} scores with respect to some objective function are also expected to be simple, hence not always optima. 
We argued that these types of optimisations are common in nature, due to the fact that many objective functions derive from generic physics and engineering-based constraints, and these constraints can be described by simple $O(1)$ equations and laws. Our result is perhaps a generalisation of the finding of Cohn and Kumar \cite{cohn2009algorithmic} that symmetrical configurations can often be built using very simple potential functions, and perhaps also the statement of Bormashenko \cite{bormashenko2022fibonacci} that the symmetry of biological structures follow the symmetry of media in which the structure is functioning. In ref.\ \cite{johnston2022symmetry} we argued that much of the symmetry in biology can be explained by the fact that simpler shapes are easier to generate via efficient genomic programs in DNA, which is a kind of `internal' argument for simplicity. Here we have shown that symmetry and simplicity can arise from a kind of `external' argument, that of the simplicity of the external constraints.

Here we also argued that given a geometry is optimal for one simple objective function, it is therefore more likely to be an optimal geometry for another simple objective function, as compared to a null expectation based on random functions. This is a rather surprising result, given that we typically conceive of exponentially large search spaces. Our argument is somewhat similar to that of Valle-Perez et al \cite{valle2018deep} who proposed that a reason for the success of deep neural networks in machine learning is that they are biased towards simpler functions, and natural functions are also tend to be simple, such that a coincidence of functions is more likely than expected by a random function null model.

The work here has focussed on discrete geometries, configurations, or sequences, such as binary strings, discrete graphs, branching patters, and lattice shapes. However, continuous curves, functions and shapes can be discretised often quite easily \cite{dingle2018input,dingle2022note}, and because of this the arguments presented here apply to a wide array of shapes and geometries. 

Looking at individual systems,  we may be able to use system specific theory to directly infer that some optimal shape for one objective function is linked to optimal shapes for some other objective function, eg employing thermodynamics to understand why high mutational robustness and thermal stability are linked in natural proteins \cite{bloom2005thermodynamic}. Importantly, such examples do not undermine our work here, in which we have attempted to show \emph{general theory} linking optima, simplicity, and coinciding optima. Whether or not particular cases can or cannot be understood with existing system specific theory, does not affect our arguments. Here we have taken an information content approach from algorithmic information theory as a theoretical framework for investigating these questions in general.\\

\noindent
{\bf Acknowledgements:} We acknowledge financial support from the Gulf University for Science and Technology for a Seed Grant (grant number 263301) and a Summer Faculty Fellowship. We thank Paris Flood and Nora Martin for discussions related to this work.\\


\bibliographystyle{unsrt}
\bibliography{OptimaRefs} 

\end{document}